\documentclass[reqno]{amsart}
\usepackage{amsmath}
\usepackage{amsfonts}
\usepackage{amssymb}
\usepackage{graphicx}
\usepackage[update,prepend]{epstopdf}

\renewenvironment{proof}[1][\proofname]{{\bfseries #1.}}{\qed}

\makeatletter
\def\specialsection{\@startsection{section}{1}
  \z@{\linespacing\@plus\linespacing}{.5\linespacing}
  {\normalfont}}
\def\section{\@startsection{section}{1}
  \z@{.7\linespacing\@plus\linespacing}{.5\linespacing}
  {\normalfont\scshape\bfseries}}
\makeatother
\newtheorem{theorem}{Theorem}
\theoremstyle{plain}

\newtheorem{corollary}[theorem]{Corollary}

\newtheorem{lemma}[theorem]{Lemma}
\newtheorem{notation}{Notation}
\newtheorem{proposition}[theorem]{Proposition}
\newtheorem{remark}{Remark}

\numberwithin{equation}{section}
\begin{document}
\title[Some Remarkable Concurrences in the Quadrilateral]{S\MakeLowercase{ome} R\MakeLowercase{emarkable} C\MakeLowercase{oncurrences} \MakeLowercase{in} \MakeLowercase{the} Q\MakeLowercase{uadrilateral}}
\author{A\MakeLowercase{ndrei} S. C\MakeLowercase{ozma}}
\date{September 28, 2012}

\begin{abstract}
We study some properties of quadrilaterals concerning concurrence of lines under few to none restrictive conditions, and obtain an  extension of a transversal theorem (see [1, page 28]) from triangles to quadrilaterals.
\end{abstract}
\maketitle

\vspace{-0.3in}

\section{Introduction} \label{intro}
\smallskip

In one of his articles (see [2]), Temistocle Birsan examines collinearity properties of a convex quadrilateral with a set of well-defined points located either on its boundary, diagonals or at the intersection of some lines. Surprisingly, given the arbitrary nature of the quadrilateral and the little hypothesis, it turns out that no less than eight lines pass through the intersection point of the diagonals.

Influenced by the aforementioned work this paper considers a more general setting, starting with a convex quadrilateral in the Euclidean plane and four points situated one on each side. Lines are drawn between these points and the vertices, and various intersection points - possibly with the diagonals - are defined. While in Birsan's article only one point's location on a line can vary and all other points are fixed with respect to it, now we change our point of view and consider a single condition for the entire system (this condition involves the ratios determined by the four points lying on the sides of the quadrilateral). 

We would expect to find a configuration poor in potentially useful properties, however it turns out to be the exact opposite: seven lines are concurrent. Afterwards we remove this condition as well and analyse the resulting configuration in the most random scenario possible.

Michael Keyton of the Illinois Mathematics and Science Academy defines in one of his articles on Euclidean geometry (see [3]) a \textit{Theorem of Mystery} as 'a result that has considerable structure with minimal hypothesis' and as we will see throughout this paper, this concept of mystery in geometry can definitely be associated with the main result. I will assume familiarity with some classical theorems in Euclidean plane geometry including but not limited to results due to Menelaus, Ceva and Van Aubel (see [4]). \\

Consider a convex quadrilateral $ABCD$ with $O$ being the intersection point of diagonals $AC$ and $BD$, and points $M,N,P,Q$ on the four sides $(AB)$, $(BC)$, $(CD)$ and $(DA)$ respectively, so that the following holds:

\begin{equation} \label{restrictive condition}
\frac{AM}{MB} \cdot \frac{BN}{NC} \cdot \frac{CP}{PD} \cdot \frac{DQ}{QA} = 1
\end{equation}

\begin{notation} 
\begin{equation}
\frac{AM}{MB} = m, \; \frac{BN}{NC} = n, \; \frac{CP}{PD} = p, \; \frac{DQ}{QA} = q \nonumber
\end{equation}
\end{notation}

\begin{proposition} 
Call the intersection points of lines $AN$ with $BQ$, and $DN$ with $CQ$, $X$ and $Z$ respectively. In the same way, lines $CM$ and $BP$ meet in $Y$, while $AP$ and $DM$ meet in $T$. Then both $XZ$ and $YT$ pass through point $O$.
\end{proposition}
\vspace{-0.1in}

\begin{figure}[h]
\begin{center}
\includegraphics[width=\linewidth]{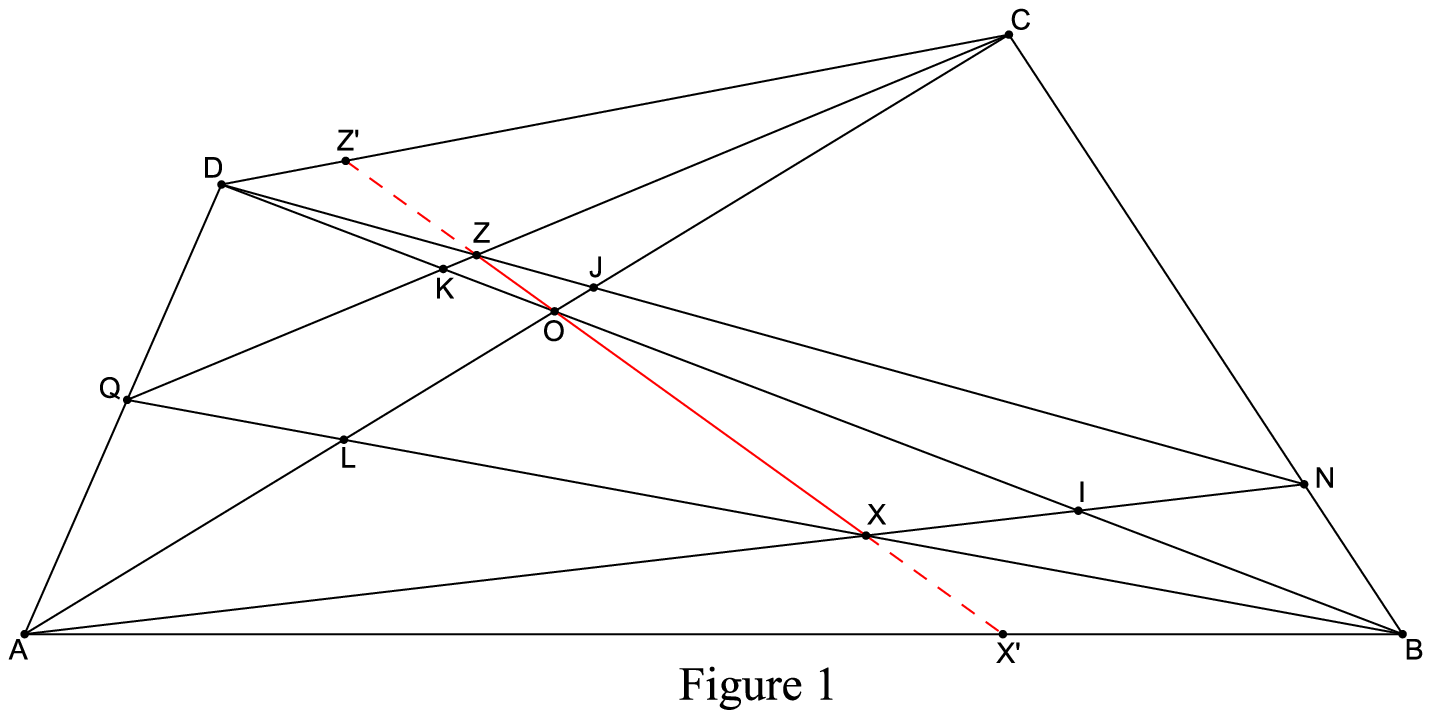}
\end{center}
\end{figure}

\noindent
\begin{proof}
\, By symmetry of the above configuration, it is enough to prove the first part of the statement, i.e. that $X$, $O$, $Z$ are collinear. Define new points:
\vspace{-0.2in}
\begin{center}
\[ {X'} = OX \cap AB, \: {Z'} = OZ \cap CD \] \vspace{-0.2in}
\[ {I} = AN \cap BD, \: {J} = DN \cap AC, \: {K} = CQ \cap BD, \: {L} = BQ \cap AC \]
\end{center}

\noindent
Th. Menelaus in $\triangle BCO$ with transversal $A-I-N$ implies: 
\[ \frac{BI}{IO} = \frac{BN}{NC} \cdot \frac{CA}{AO} = n \cdot \frac{CA}{AO} \]
Similar applications of Th. Menelaus in triangles $\triangle BCO$ with transversal $D-J-N$ and $\triangle ADO$ with transversals $B-L-Q$ and $C-K-Q$ give:
\[ \frac{CJ}{JO} = \frac{1}{n} \cdot \frac{BD}{DO}; \; \frac{AL}{LO} = \frac{1}{q} \cdot \frac{BD}{BO}; \; \frac{DK}{KO} = q \cdot \frac{AC}{CO}; \]
We now apply Th. Ceva twice in triangles $\triangle AOB$ and $\triangle COD$ to get:
\[ \frac{AX'}{X'B} = \frac{AL}{LO} \cdot \frac{OI}{IB} = \frac{1}{nq} \cdot \frac{BD}{BO} \cdot \frac{AO}{AC}, \; \text{ and}\]
\[ \frac{CZ'}{Z'D} = \frac{CJ}{JO} \cdot \frac{OK}{KD} = \frac{1}{nq} \cdot \frac{BD}{DO} \cdot \frac{CO}{AC}. \]
Combining the last two equations leads to:
\begin{equation} \frac{CZ'}{Z'D} = \frac{AX'}{X'B} \cdot \frac{BO}{OD} \cdot \frac{CO}{OA}. \end{equation}
\smallskip

In order to achieve the collinearity of the three points we define ${Z"} = XO \cap CD$ and only need prove that it coincides with the point $Z'$. By using the equality of the angle measures around point $O$ and the Sine Law we easily reach the conclusion:
\[ \frac{\sin \widehat{AOX'}}{\sin \widehat{BOX'}} = \frac{\sin \widehat{COZ"}}{\sin \widehat{DOZ"}} \Rightarrow \frac{AX'}{AO} \cdot \frac{BO}{BX'} = \frac{CZ"}{CO} \cdot \frac{DO}{DZ"} \Rightarrow \frac{CZ"}{DZ"} = \frac{AX'}{BX'} \cdot \frac{BO}{OA} \cdot \frac{CO}{OD} = \frac{CZ'}{DZ'}. \]
\end{proof}

\section{The Point of the Seven Lines}

Consider the configuration established at the beginning of the previous section, and define new points: 
\[ {A'} = BP \cap DN, \, {B'} = AP \cap CQ, \, {C'} = BQ \cap DM, \, {D'} = AN \cap CM, \]
\[ {F_1} = BD' \cap AC, \, {G_1} = CA' \cap BD, \, {F_2} = DB' \cap AC, \, {G_2} = AC' \cap BD. \]
\vspace{-0.2in}

\begin{figure}[h]
\begin{center}
\includegraphics[width=\linewidth]{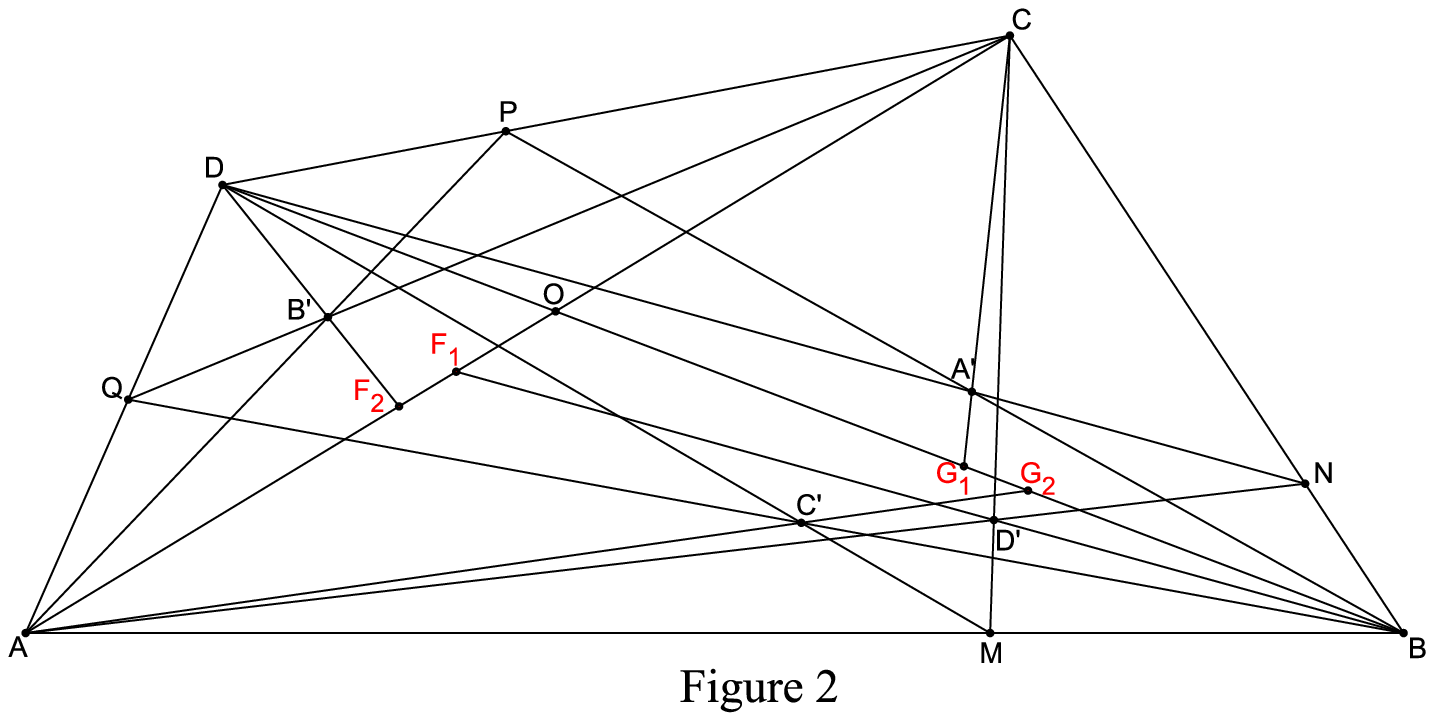}
\end{center}
\end{figure}

\noindent
Using Th. Ceva in triangles $\triangle ABC$ and $\triangle ADC$ we find that:
\[ \frac{AF_{1}}{F_{1}C} = \frac{AM}{MB} \cdot \frac{BN}{NC} = m \cdot n = \frac{m \cdot n \cdot p \cdot q}{p \cdot q} = \frac{1}{p \cdot q} = \frac{DP}{PC} \cdot \frac{AQ}{QD} = \frac{AF_{2}}{F_{2}C}. \] 
\noindent
$\therefore$ Therefore the two points coincide and $F_{1} \equiv F_{2} \equiv F$. Similarly, $G_{1} \equiv G_{2} \equiv G$.
\bigskip

In order to state and prove the main result of this section, we first need to cover two auxiliary lemmas:

\begin{lemma}
Lines $AA'$, $BB'$, $CC'$ and $DD'$ intersect the segment $[FG]$ in the same point, call it ${E}$.
\end{lemma}
\vspace{-0.2in}

\begin{figure}[h]
\begin{center}
\includegraphics[width=\linewidth]{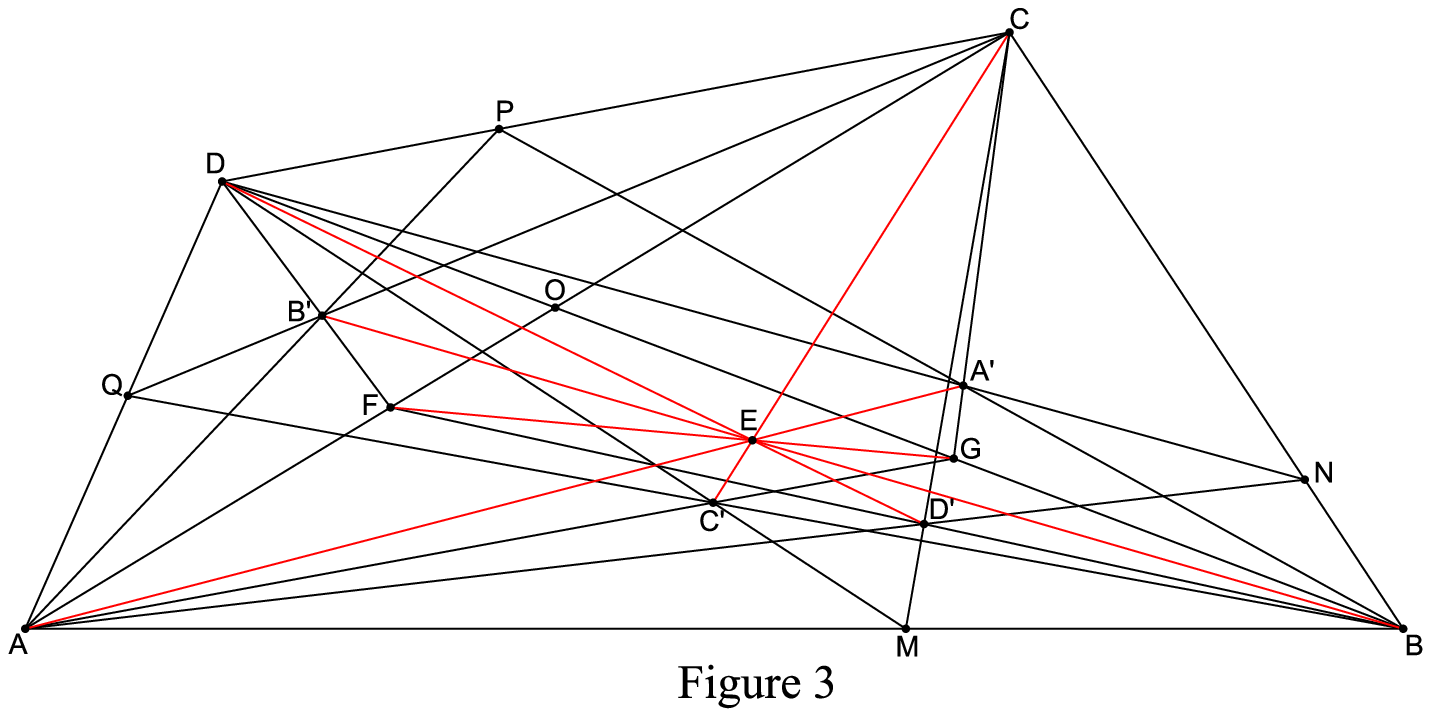}
\end{center}
\end{figure}

\noindent
\begin{proof} 
\, First of all, we notice that since $A'$ lies on the open segment $(CG)$ and $F$ on $(AC)$, then the two segments $[FG]$ and $[AA']$ do in fact intersect. The same argument holds for $[BB'], [CC'], [DD']$. \\

\noindent
\textit{Step 1.} Lines $CC', DD', FG$ are concurrent. \\

\noindent
Firstly, let ${E} = CC' \cap FG$ and ${E'} = DD' \cap FG$. It is enough to prove that the two points coincide. Using Th. Menelaus and Th. Ceva in several triangles leads to
\begin{equation}
\frac{AC'}{C'G} = \frac{AM}{MB} \cdot \frac{BD}{DG} = \frac{AM}{MB} \cdot \left[1 + \frac{BG}{GD}\right] = \frac{AM}{MB} \cdot \left[1 + \frac{BN}{NC} \cdot \frac{CP}{PD} \right] = m (1 + np), \nonumber
\end{equation}
\begin{equation}
\frac{CF}{FA} = \frac{CN}{NB} \cdot \frac{BM}{MA} = \frac{1}{mn} \Rightarrow \frac{CF}{AC} = \frac{1}{mn + 1}, \text{ so } \frac{FE}{EG} = \frac{FC}{AC} \cdot \frac{AC'}{C'G} = \frac{m(1+np)}{1+mn}. \nonumber
\end{equation} 

\noindent
Next we consider the ratio:
\begin{equation} \label{F-E-G ratio}
\frac{FE'}{E'G} = \frac{FD'}{D'B} \cdot \frac{BD}{DG} = \frac{FC}{CA} \cdot \frac{AM}{MB} \cdot \frac{BD}{DG} = \frac{1}{1+mn} \cdot m \cdot (1+np) = \frac{m(1+np)}{1+mn}
\end{equation}
therefore $\frac{FE'}{E'G} = \frac{FE}{EG}$, points $E$ and $E'$ are one and the same and concurrence follows. \\

\noindent
\textit{Step 2.} Lines $AA', CC', FG$ are concurrent. \\

\noindent
We use Th. Menelaus and Th. Ceva repeatedly as follows:
\[ \frac{CA'}{A'G} \cdot \frac{GC'}{C'A} \cdot \frac{AF}{FC} = \left[\frac{CN}{NB} \cdot \frac{BD}{DG}\right] \cdot \left[ \frac{BG}{BD} \cdot \frac{DQ}{QA}\right] \cdot \left[\frac{AQ}{QD} \cdot \frac{DP}{PC} \right] = \frac{1}{n} \cdot \frac{BG}{GD} \cdot \frac{1}{p} = 1, \]
hence the reciprocal of Th. Ceva implies the concurrence of the three lines, which are in fact cevians in $\triangle ACG$. \\

\noindent
\textit{Step 3.} Lines $BB', DD', FG$ are concurrent. \\

\noindent
Using the same method as above, we obtain:
\[ \frac{DB'}{B'F} \cdot \frac{FD'}{D'B} \cdot \frac{BG}{GD} = \left[\frac{DP}{PC} \cdot \frac{CA}{AF}\right] \cdot \left[ \frac{FA}{AC} \cdot \frac{CN}{NB}\right] \cdot \left[\frac{BN}{NC} \cdot \frac{CP}{PD} \right] = 1, \]
hence the reciprocal of Th. Ceva implies the concurrence of the three lines, which are cevians in $\triangle BDF$. \\

\noindent
$\therefore$ In conclusion, lines $AA', BB', CC', DD' \text{ and } FG$ are concurrent in $E$.
\end{proof}
\bigskip

\begin{lemma}
Lines $MP \text{ and } NQ$ intersect the segment $FG$ in $E$.
\end{lemma}

\noindent
\begin{proof}
\, We already know from the previous lemma that $E$ is located at the intersection of lines $AA', BB', CC', DD' \text{ and } FG$. Using Th. Menelaus in triangles $\triangle BCP$ and $\triangle ADP$ with transversals $D-A'-N$ and $C-B'-Q$ respectively,
\[ \frac{BA'}{A'P} \cdot \frac{PB'}{B'A} \cdot \frac{AM}{MB} = \left[ \frac{BN}{NC} \cdot \frac{CD}{DP} \right] \cdot \left[ \frac{PC}{CD} \cdot \frac{DQ}{QA} \right] \cdot m = n(1+p) \cdot \frac{pq}{1+p} \cdot m = m \cdot n \cdot p \cdot q = 1, \]
hence by the reciprocal of Th. Ceva, $AA', BB' \text{ and } MP$ are cevians in $\triangle ABP$. \\

\noindent
$\therefore$ In conclusion, point $E$ lies on segment $[MP]$. Working in a similar fashion we can prove that $BB', CC' \text{ and } NQ$ are cevians in $\triangle BCQ$, so $E$ lies on segment $[NQ]$ as well, and the conclusion follows.
\end{proof}
\bigskip

\begin{theorem} 
[Concurrence of seven lines] The lines $AA', BB', CC', DD', MP, NQ$ and $FG$ are concurrent.
\end{theorem}

\noindent
\begin{proof}
\, The previous two results assure us that the seven lines in the statement intersect in $E$.
\end{proof}

\begin{figure}[h]
\begin{center}
\includegraphics[width=\linewidth]{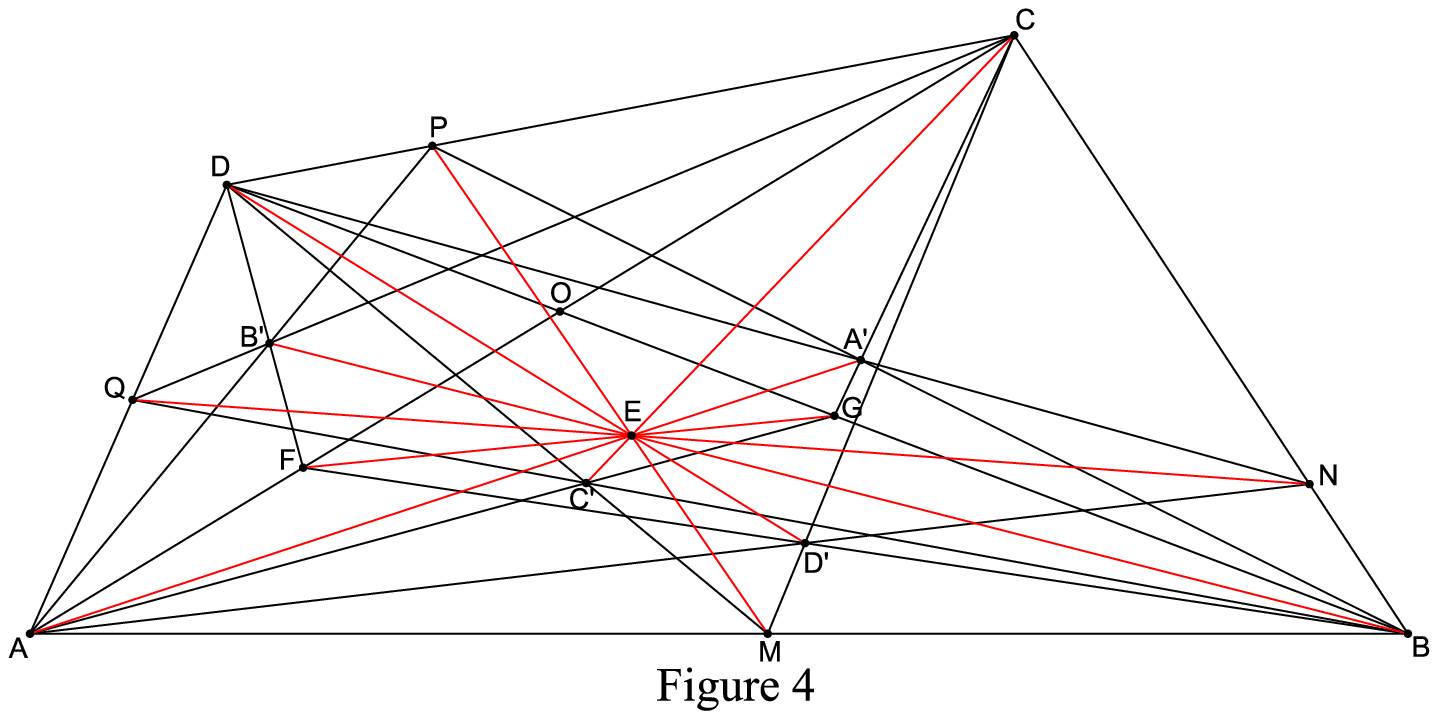}
\end{center}
\end{figure}

\begin{theorem}
[Extension of a Transversal Theorem] Consider a convex quadrilateral $ABCD$ and points $M, N, P, Q$ on sides $(AB), (BC), (CD), (DA)$ respectively, satisfying \, $\frac{AM}{MB} \cdot \frac{BN}{NC} \cdot \frac{CP}{PD} \cdot \frac{DQ}{QA} = 1.$ Define $E = MP \cap NQ$, then:

\[ (1) \, \frac{ME}{EP} = \frac{AQ}{QD} \cdot \frac{MB}{AB} + \frac{BN}{NC} \cdot \frac{MA}{AB}; \hspace{0.3in} (2) \, \frac{NE}{EQ} = \frac{BM}{MA} \cdot \frac{NC}{BC} + \frac{CP}{PD} \cdot \frac{NB}{BC}. \]
\end{theorem}
\vspace{-0.2in}

\begin{figure}[hb]
\begin{center}
\includegraphics[width=\linewidth]{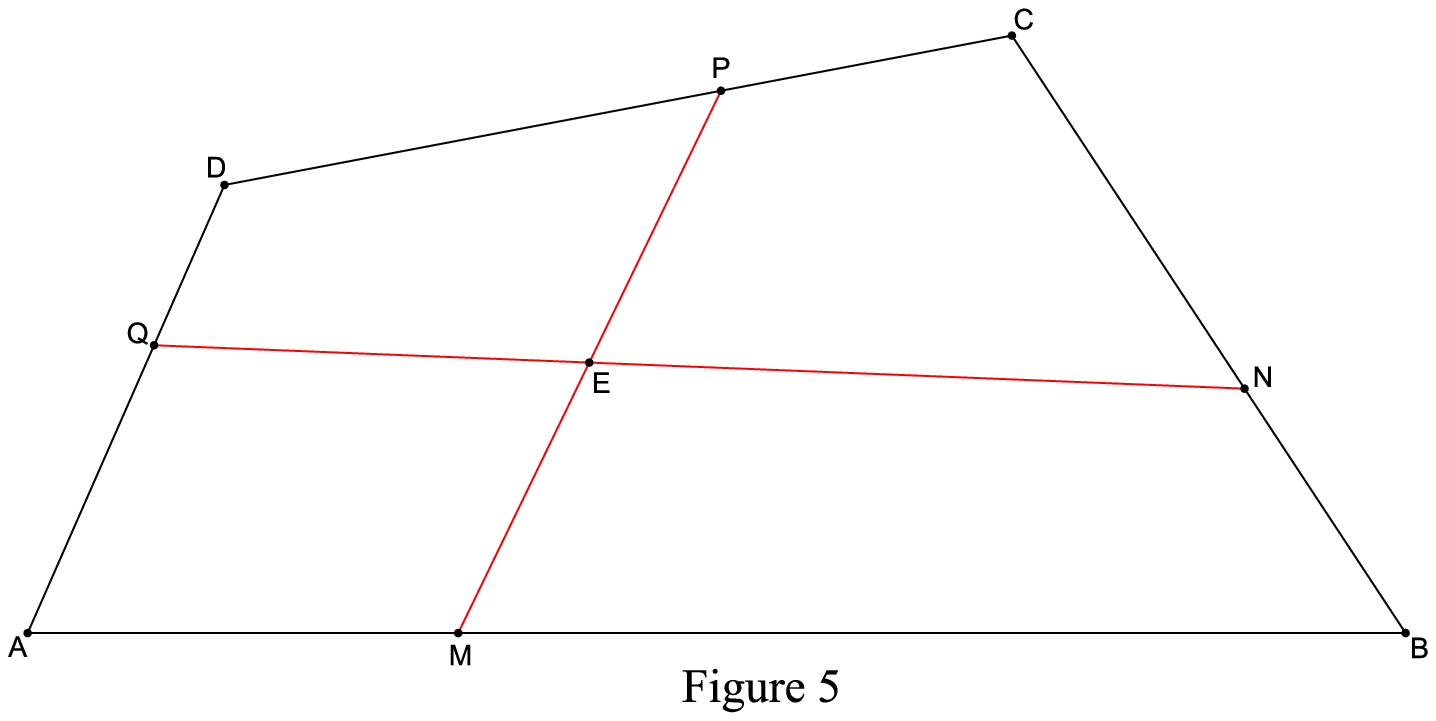}
\end{center}
\end{figure}

\noindent
\begin{proof}
\, Using Th. Van Aubel in triangle $\triangle ABP$ and some results derived in the proof of \textit{Lemma 3}, we find:
\[ \frac{PE}{EM} = \frac{PB'}{B'A} + \frac{PA'}{A'B}  = \frac{pq}{p+1} + \frac{1}{n(p+1)} = \frac{m(npq+1)}{mn(p+1)} = \frac{m+1}{mn(p+1)}. \]

\[ \text{In other words, } \frac{ME}{EP} = \frac{mn(p+1)}{m+1} = \frac{mnp}{m+1} + \frac{mn}{m+1} = \frac{1}{q} \cdot \frac{1}{m+1} + \frac{m}{m+1} \cdot n \]
\[ \Rightarrow \frac{ME}{EP} = \frac{AQ}{QD} \cdot \frac{MB}{AB} + \frac{BN}{NC} \cdot \frac{MA}{AB}, \text{ as required.} \]

\[ \text{Similarly, } \frac{NE}{EQ} = \frac{np(q+1)}{n+1} = \frac{1}{m} \cdot \frac{1}{n+1} + p \cdot \frac{n}{n+1} = \frac{BM}{MA} \cdot \frac{NC}{BC} + \frac{CP}{PD} \cdot \frac{NB}{BC}. \]
\end{proof}
\bigskip

\begin{corollary}
If the points $M$, $N$, $P$ and $Q$ are located on the four sides of the quadrilateral so that $\frac{AQ}{QD} = \frac{BN}{NC}$ and $\frac{AM}{MB} = \frac{DP}{PC}$, then condition (\ref{restrictive condition}) is satisfied and we find that $\frac{ME}{EP} = \frac{AQ}{QD}$, $\frac{NE}{EQ} = \frac{BM}{MA}$. In particular, if $M$, $N$, $P$, $Q$ are midpoints, then $ME = EP$ and $NE = EQ$.
\end{corollary}
\medskip

Next we examine two applications of \textit{Theorem 4} by choosing particular locations of points $F$ and/or $G$ on diagonals $AC$ and $BD$ respectively.
\bigskip

\begin{corollary}
Let $ABCD$ be a convex quadrilateral and $M$, $N$, $P$, $Q$ on the four sides, $(AB)$, $(BC)$, $(CD)$ and $(DA)$ respectively, so that \, $\frac{AM}{MB} \cdot \frac{BN}{NC} \cdot \frac{CP}{PD} \cdot \frac{DQ}{QA} = 1.$ Then $MP$, $NQ$ and $AC$ are concurrent if and only if $DM$, $BQ$ and $AC$ are.
\end{corollary}

\begin{figure}[h]
\begin{center}
\includegraphics[width=\linewidth]{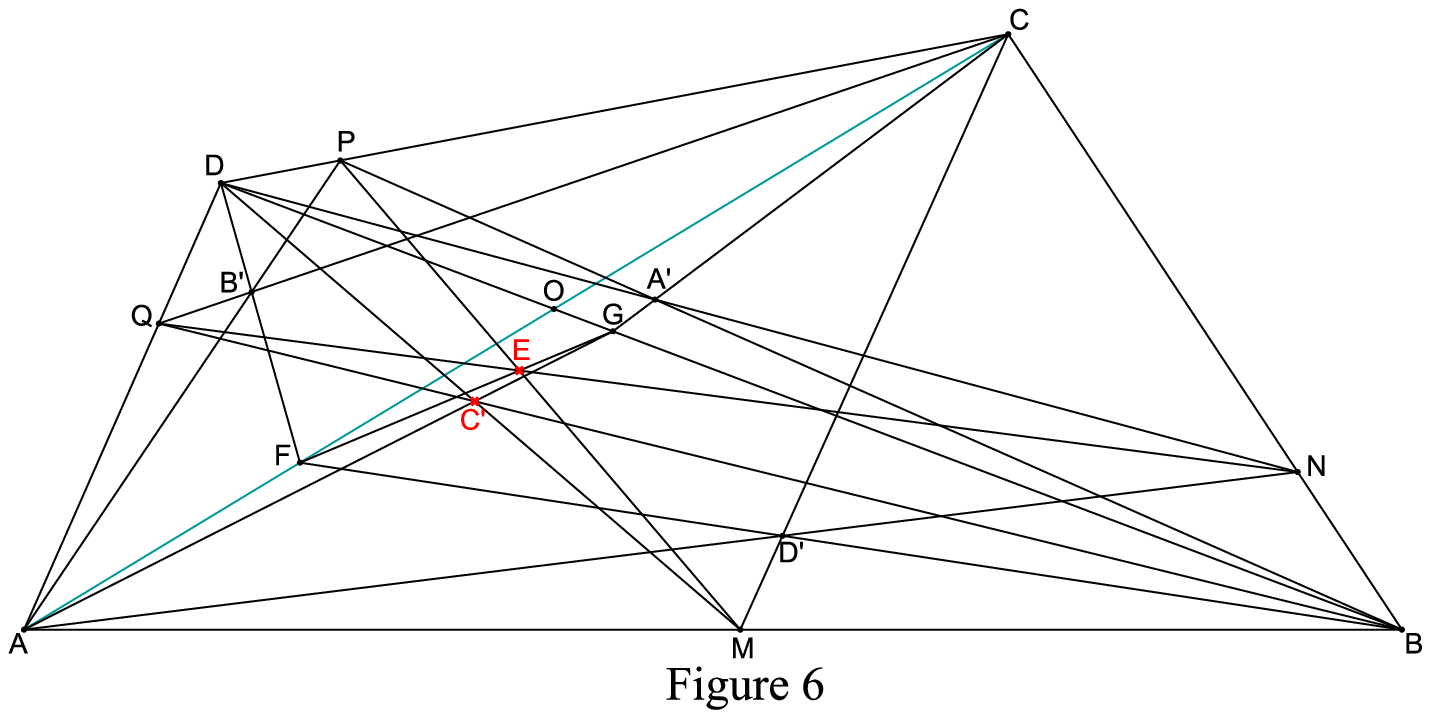}
\end{center}
\end{figure}

\noindent
\begin{proof}
\, We know that $E \in (CC')$ so $E \in (AC)$ if and only if $C' \in (AC)$. Therefore, taking into account that ${E} = MP \cap NQ$, we have that $MP \cap NQ \cap AC \neq \emptyset$ if and only if $E \in (AC)$ if and only if $C' \in (AC)$ if and only if $BQ \cap DM \cap AC \neq \emptyset$, because $BQ \cap DM = C'$.
\end{proof}
\bigskip

\begin{corollary}
Consider a convex quadrilateral $ABCD$ and points $M$ and $N$ on sides $(AB)$ and $(BC)$ respectively. Define new points: ${S} = DM \cap AC$, ${T} = DN \cap AC$, ${Q} = BS \cap AD$, ${P} = BT \cap CD$. Then lines $MP$, $NQ$ and $AC$ are concurrent.
\end{corollary}

\vspace{-0.2in}

\begin{figure}[h]
\begin{center}
\includegraphics[width=\linewidth]{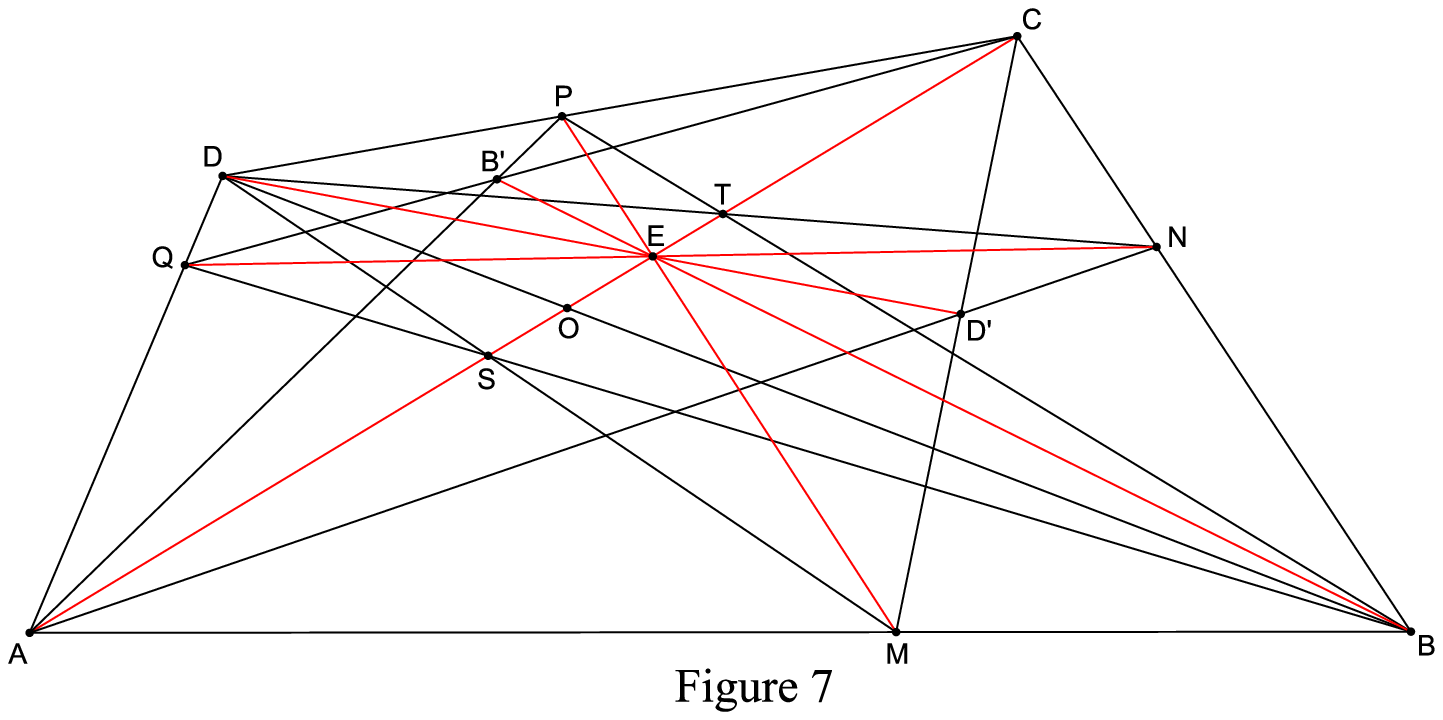}
\end{center}
\end{figure}

\noindent
\begin{proof}
\, Follows immediately from the previous corollary by noticing that our construction implies $\frac{AM}{MB} \cdot \frac{BN}{NC} \cdot \frac{CP}{PD} \cdot \frac{DQ}{QA} = 1$ (use Ceva twice in triangles $\triangle ABD$ and $\triangle BCD$). In fact, $BB'$ and $DD'$ pass through the same intersection point.
\end{proof} 
\vspace{0.3cm}

\section{Concurrence in a General Setting}

Consider a convex quadrilateral $ABCD$ and four points $M, N, P, Q$ on the sides $(AB)$, $(BC)$, $(CD)$, $(DA)$ respectively, chosen at random. Define points: ${O} = AC \cap BD$, ${A'} = BP \cap DN$, ${B'} = AP \cap CQ$, ${C'} = BQ \cap DM$, ${D'} = AN \cap CM$, ${F_1} = BD' \cap AC$, ${G_1} = CA' \cap BD$, ${F_2} = DB' \cap AC$, ${G_2} = AC' \cap BD$.
\medskip

\begin{notation}
\begin{equation} \label{non-restrictive condition}
\frac{AM}{MB} \cdot \frac{BN}{NC} \cdot \frac{CP}{PD} \cdot \frac{DQ}{QA} = m \cdot n \cdot p \cdot q = \gamma
\end{equation}
\end{notation}

\noindent
Using Th. Ceva repeatedly in four triangles, we find that:
\[ \frac{AF_{1}}{F_{1}C} = \frac{AM}{MB} \cdot \frac{BN}{NC} = mn; \, \frac{BG_{1}}{G_{1}D} = \frac{BN}{NC} \cdot \frac{CP}{PD} = np; \, \frac{AF_{2}}{F_{2}C} = \frac{AQ}{QD} \cdot \frac{DP}{PC} = \frac{1}{pq} = \frac{mn}{\gamma}; \]
\[ \frac{BG_{2}}{G_{2}D} = \frac{BM}{MA} \cdot \frac{AQ}{QD} = \frac{1}{mq} = \frac{np}{\gamma}; \]
\begin{equation}
\Rightarrow \frac{AF_{1}}{F_{1}C} \cdot \frac{CF_{2}}{F_{2}A} = mn \cdot \frac{\gamma}{mn} = \gamma = np \cdot \frac{\gamma}{np} = \frac{BG_{1}}{G_{1}D} \cdot \frac{DG_{2}}{G_{2}B}. 
\end{equation}
\smallskip

\noindent
Thus we encounter three possible ranges of values for $\gamma$:
\begin{itemize}
\item \textit{case 1}. $\gamma = 1 \Rightarrow F_{1} \equiv F_{2} \text{ and } G_{1} \equiv G_{2}$, i.e. the \textit{Section 2} setting.
\item \textit{case 2}. $\gamma < 1 \Rightarrow A, F_{1}, F_{2}, C \text{ and } B, G_{1}, G_{2}, D$ appear on their respective diagonals in this order.
\item \textit{case 3}. $\gamma > 1 \Rightarrow A, F_{2}, F_{1}, C \text{ and } B, G_{2}, G_{1}, D$ appear on their respective diagonals in this order.
\end{itemize}
\medskip

Before stating and proving the main result of this section, we first need two auxiliary lemmas:

\begin{lemma}
The triplets of lines $\{DD', AA', F_{1}G_{1}\}$, $\{AA', BB', G_{1}F_{2}\}$, $\{BB', CC', F_{2}G_{2}\}$ and $\{CC', DD', G_{2}F_{1}\}$ meet in points $M_{1}$, $N_{1}$, $P_{1}$ and $Q_{1}$, respectively.
\end{lemma}
\vspace{-0.1in}

\begin{figure}[h]
\begin{center}
\includegraphics[width=\linewidth]{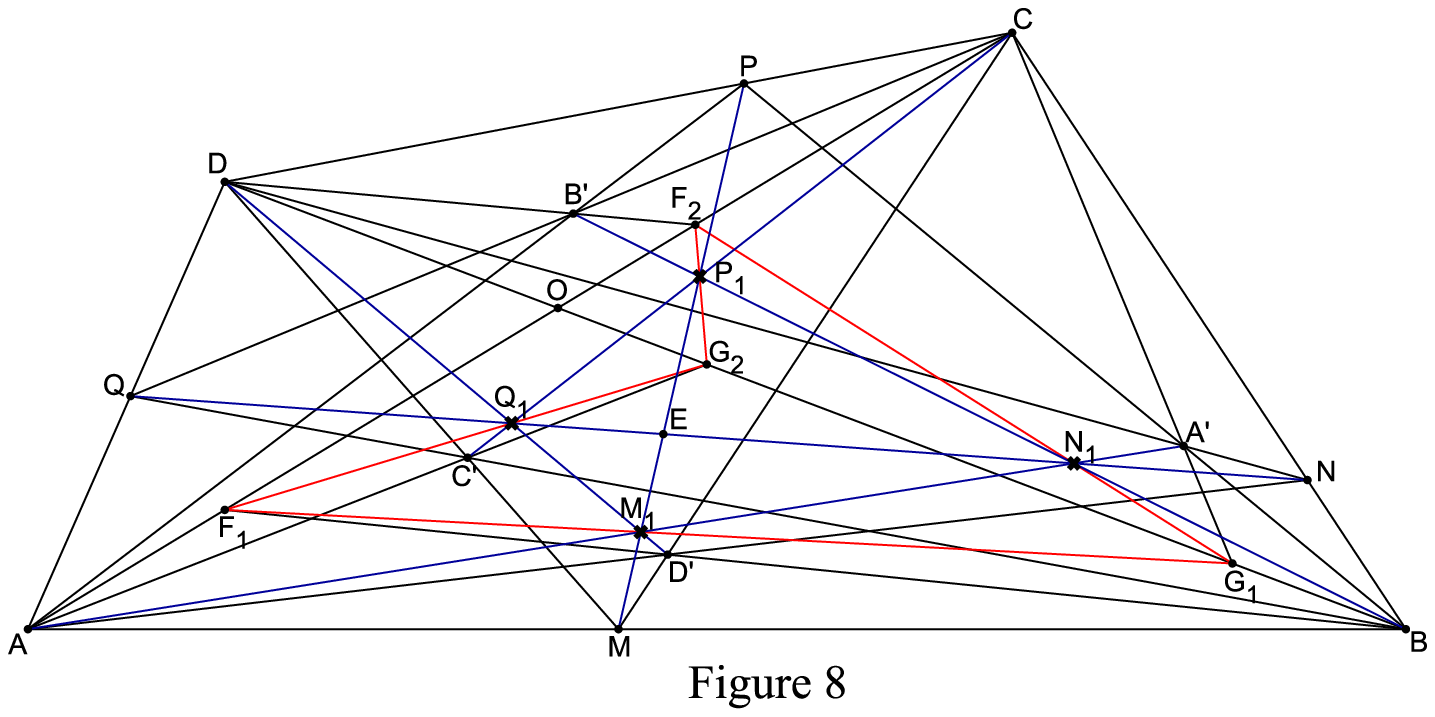}
\end{center}
\end{figure}

\noindent
\begin{proof}
\, By symmetry of the configuration, it is enough to show that $AA' \cap DD' \cap F_{1}G_{1} \neq \emptyset$, and the remaining concurrences follow. Set $M_{1}$ to be the intersection point of $F_{1}G_{1}$ with $AA'$, and $M_{1}'$ to be the intersection point of $F_{1}G_{1}$ with $DD'$.

\[ \text{Using Th. Menelaus in } \triangle F_{1}G_{1}C \text{ with transversal } A-M_{1}-A': \frac{F_{1}M_{1}}{M_{1}G_{1}} = \frac{F_{1}A}{AC} \cdot \frac{CA'}{A'G_{1}}; \]
\[ \text{Th. Ceva in } \triangle BCA: \frac{AF_{1}}{F_{1}C} = \frac{AM}{MB} \cdot \frac{BN}{NC} = mn \Rightarrow \frac{AF_{1}}{AC} = \frac{mn}{mn+1}; \]
\[ \text{Th. Van Aubel in } \triangle BCD: \frac{CA'}{A'G_{1}} = \frac{CN}{NB} + \frac{CP}{PD} = \frac{1}{n} + p = \frac{np+1}{n}; \]
\begin{equation}
\Rightarrow \frac{F_{1}M_{1}}{M_{1}G_{1}} = \frac{mn}{mn+1} \cdot \frac{np+1}{n} = \frac{m(np+1)}{mn+1}
\end{equation}

\[ \text{Using Th. Menelaus in } \triangle F_{1}G_{1}B \text{ with transversal } D-M_{1}'-D': \frac{F_{1}M_{1}'}{M_{1}'G_{1}} = \frac{F_{1}D'}{D'B} \cdot \frac{BD}{DG_{1}}; \]
\[ \text{Th. Ceva in } \triangle BCD: \frac{BG_{1}}{G_{1}D} = \frac{BN}{NC} \cdot \frac{CP}{PD} = np \Rightarrow \frac{BD}{G_{1}D} = np + 1; \]
\[ \text{Th. Van Aubel in } \triangle ABC: \frac{BD'}{D'F_{1}} = \frac{BM}{MA} + \frac{BN}{NC} = \frac{1}{m} + n = \frac{mn+1}{m}; \]
\begin{equation}
\Rightarrow \frac{F_{1}M_{1}'}{M_{1}'G_{1}} = \frac{m}{mn+1} \cdot (np+1) = \frac{m(np+1)}{mn+1} = \frac{F_{1}M_{1}}{M_{1}G_{1}}
\end{equation}
\smallskip

\noindent
$\therefore$ In conclusion, $M_{1}' \equiv M_{1}$ and hence $F_{1}G_{1} \cap AA' \cap DD' = M_{1}$ as required. 
\end{proof}
\bigskip

\begin{lemma}
If points $M_{1}$, $N_{1}$, $P_{1}$ and $Q_{1}$ are as defined above, then the lines $MP$ and $NQ$ intersect $\{F_{1}G_{1}, F_{2}G_{2}\}$ and $\{G_{1}F_{2}, G_{2}F_{1}\}$ in $\{M_{1}, P_{1}\}$ and $\{N_{1}, Q_{1}\}$, respectively.
\end{lemma}

\noindent
\begin{proof}
\, Again it is easy to see that proving the concurrence of $MP$ and $F_{1}G_{1}$ in $M_{1}$ is enough, as the rest follows from symmetry. Define a new point, $M' = PM_{1} \cap AB$, and we only need prove that $M'$ and $M$ coincide to reach the desired conclusion. \\

\noindent
Using Th. Menelaus in $\triangle CAA'$ with transversal $F_{1}-M_{1}-G_{1}$:
\[ \frac{AM_{1}}{M_{1}A'} = \frac{AF_{1}}{F_{1}C} \cdot \frac{CG_{1}}{G_{1}A'} = mn \cdot \frac{np+n+1}{n} = m(np+n+1); \]
\[ \frac{A'B}{A'P} = \frac{BN}{NC} \cdot \frac{CD}{DP} = n(p+1) \Rightarrow \frac{A'P}{PB} = \frac{1}{np+n+1}; \]
Th. Menelaus in $\triangle BAA'$ with transversal $P-M_{1}-M'$:
\[ \frac{AM'}{M'B} = \frac{AM_{1}}{M_{1}A'} \cdot \frac{A'P}{PB} = m(np+n+1) \cdot \frac{1}{np+n+1} = m = \frac{AM}{MB}. \]

\noindent
$\therefore$ Therefore, $M' \equiv M$ so $M_{1} \in (MP)$ as required.
\end{proof}
\bigskip

\begin{theorem}
The four lines $F_{1}G_{1}, DD', AA', MP$ are concurrent. Similarly, $G_{1}F_{2}, AA', BB', NQ$ are concurrent, $F_{2}G_{2}, BB', CC', MP$ are concurrent, and finally $G_{2}F_{1}, CC', DD', NQ$ are concurrent.
\end{theorem}

\noindent
\begin{proof}
\, The previous two lemmas assure us that the concurrences in the statement do hold.
\end{proof}
\bigskip

\begin{corollary}
If $\gamma \, (= mnpq) = 1$, then $F_{1} \equiv F_{2} \equiv F$ and $G_{1} \equiv G_{2} \equiv G$ (see beginning of \textit{Section 3}). This implies that the segments $F_{1}G_{1}$, $F_{1}G_{2}$, $F_{2}G_{1}$, $F_{2}G_{2}$ are the same with $FG$, and the previous theorem assures us that $M_{1} \equiv N_{1} \equiv P_{1} \equiv Q_{1} \equiv E$, because every two points are located on two common lines so they must coincide. Therefore the lines $AA'$, $BB'$, $CC'$, $DD'$, $MP$, $NQ$ and $FG$ are concurrent, and \textit{Theorem 4} is a particular case of \textit{Theorem 11}.
\end{corollary}
\bigskip

\noindent
\textbf{Observations.} \\

\noindent
$(1)$ The symmetric configuration of the system allows us to compute other ratios immediately, by permuting the vertices around the quadrilateral $(\text{eg. } A \mapsto B \mapsto C \mapsto D \mapsto A)$ and the positions of $M$, $N$, $P$, $Q$ on the four sides of the quadrilateral $(m \mapsto n \mapsto p \mapsto q \mapsto m)$ at the same time. \\

\noindent
$(2)$ We already found that $\frac{F_{1}M_{1}}{M_{1}G_{1}} = \frac{m(np+1)}{mn+1}$, and by previous observation we obtain other ratios: $\frac{G_{1}N_{1}}{N_{1}F_{2}} = \frac{n(pq+1)}{np+1}$ etc. therefore:
\[ \frac{F_{1}M_{1}}{M_{1}G_{1}} \cdot \frac{G_{1}N_{1}}{N_{1}F_{2}} \cdot \frac{F_{2}P_{1}}{P_{1}G_{2}} \cdot \frac{G_{2}Q_{1}}{Q_{1}F_{1}} = \prod_{cyc.} \frac{m(np+1)}{mn+1} = m \cdot n \cdot p \cdot q = \frac{AM}{MB} \cdot \frac{BN}{NC} \cdot \frac{CP}{PD} \cdot \frac{DQ}{QA}. \] \\

\noindent
$(3)$ Using Th. Menelaus in triangle $\triangle MPB$ with transversal $A-M_{1}-A'$:
\[ \frac{MM_{1}}{M_{1}P} = \frac{MA}{AB} \cdot \frac{BA'}{A'P} = \frac{MA}{AB} \cdot \frac{BN}{NC} \cdot \frac{CD}{DP} = \frac{m}{m+1} \cdot n \cdot (p+1) = \frac{mn(p+1)}{m+1}, \]
\[ \Rightarrow \frac{MM_{1}}{M_{1}P} = \frac{mnp}{m+1} + \frac{mn}{m+1} = \gamma \cdot \frac{1}{q} \cdot \frac{1}{m+1} + n \cdot \frac{m}{m+1}, \]
\begin{equation}
\Rightarrow \frac{MM_{1}}{M_{1}P} = \gamma \cdot \frac{AQ}{QD} \cdot \frac{MB}{AB} + \frac{BN}{NC} \cdot \frac{MA}{AB}. 
\end{equation}
\smallskip

\noindent
We can determine the other ratios: $\frac{NN_{1}}{N_{1}Q}$, $\frac{PP_{1}}{P_{1}M}$ and $\frac{QQ_{1}}{Q_{1}N}$ by using \textit{Observation (1)}.
\bigskip

\begin{lemma}
Consider a convex quadrilateral $ABCD$ and points $M$, $N$, $P$, $Q$ on the four sides, $(AB)$, $(BC)$, $(CD)$ and $(DA)$ respectively. Let $NQ$ intersect $AB$ in ${R}$ and set $\frac{RA}{RB} = r$ (if $NQ \parallel AB$ define $r$ to be $1$). Then, \[ \frac{ME}{EP} = \frac{mn(p+1)}{m+1} \cdot \frac{r+m}{\gamma r + m}. \]
\end{lemma}
\smallskip

\noindent
\begin{proof}
\, First of all define new points: $X = DM \cap NQ$, $Y = CM \cap NQ$, and call $\frac{DX}{XM} = x$, $\frac{CY}{YM} = y$. We encounter three cases, depending on the location of $R$. \\

\noindent
$\bullet$ \textit{case 1.} $NQ \cap (BA = R$
\vspace{-0.2in}
\begin{figure}[h]
\begin{center}
\includegraphics[width=\linewidth]{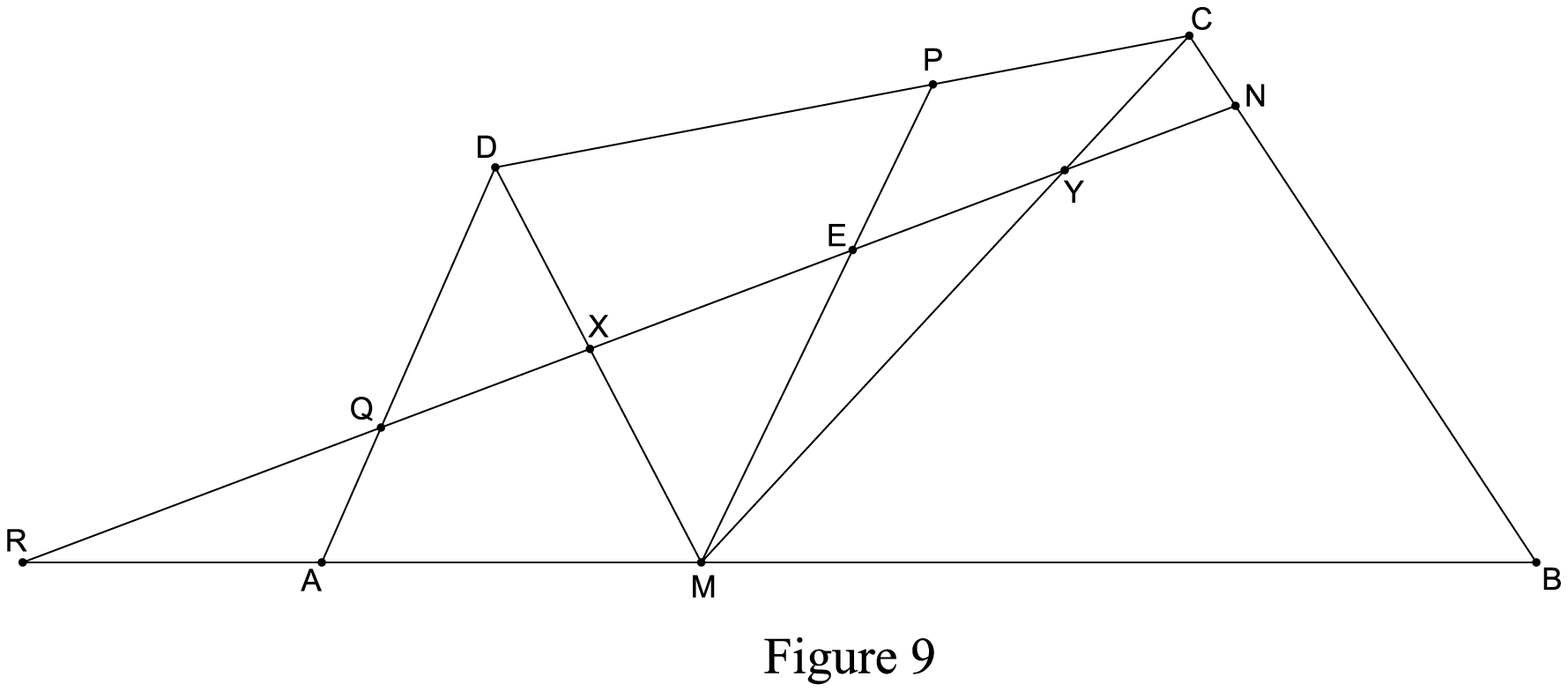}
\end{center}
\end{figure}
\vspace{-0.2in}

\noindent
Using the Transversal Theorem (see $[1]$) in triangle $\triangle CDM$, we find that:
\[ \frac{PE}{EM} = \frac{DX}{XM} \cdot \frac{CP}{CD} + \frac{CY}{YM} \cdot \frac{DP}{CD} = x \cdot \frac{p}{p+1} + y \cdot \frac{1}{p+1}; \]

\noindent
Th. Menelaus in triangle $\triangle ADM$ with transversal $R-Q-X$: 
\[ x = \frac{DX}{XM} = \frac{DQ}{QA} \cdot \frac{AR}{RM} = q \cdot \frac{AR}{RM} = q \cdot \frac{AR}{AR + AM}, \]
\[ \Rightarrow x = \frac{q \cdot AR}{AR + \frac{m}{m+1} \cdot AB} = \frac{q \cdot AR}{AR + \frac{m}{m+1} \cdot \frac{1-r}{r} \cdot AR} = \frac{q}{1 + \frac{m}{m+1} \cdot \frac{1-r}{r}} = \frac{qr(m+1)}{r+m}. \]

\noindent
On the other hand, Th. Menelaus in triangle $\triangle BCM$ with transversal $R-Y-N$:
\[ y = \frac{CY}{YM} = \frac{CN}{NB} \cdot \frac{BR}{RM} = \frac{1}{n} \cdot \frac{RB}{RA} \cdot \frac{RA}{RM} = \frac{1}{rn} \cdot \frac{r(m+1)}{r+m} = \frac{m+1}{n(r+m)}. \]

\noindent
\[ \text{Therefore, } \frac{PE}{EM} = \frac{qr(m+1)}{r+m} \cdot \frac{p}{p+1} + \frac{m+1}{n(r+m)} \cdot \frac{1}{p+1} = \frac{(m+1)(npqr+1)}{n(p+1)(r+m)}, \]
\[ \Rightarrow \frac{ME}{EP} = \frac{n(p+1)(r+m)}{(m+1)(npqr+1)} = \frac{mn(p+1)(r+m)}{(m+1)(m + \gamma r)} = \frac{mn(p+1)}{m+1} \cdot \frac{r+m}{\gamma r + m}. \]
\smallskip

\noindent
$\bullet$ \textit{case 2.} $NQ \cap (AB = R$
\vspace{-0.2in}

\begin{figure}[h]
\begin{center}
\includegraphics[width=\linewidth]{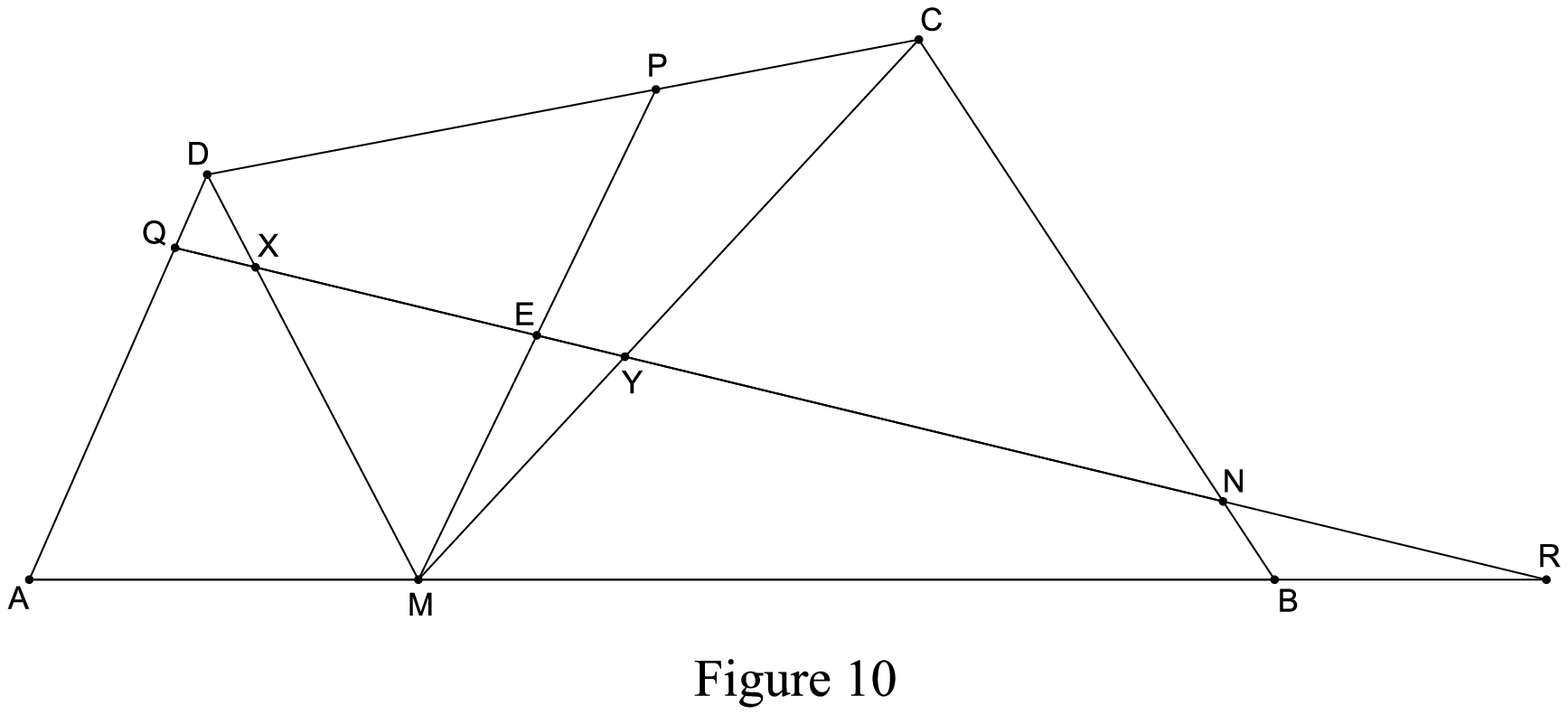}
\end{center}
\end{figure}
\vspace{-0.2in}

\noindent
\[ \text{We already know that } \frac{PE}{EM} = x \cdot \frac{p}{p+1} + y \cdot \frac{1}{p+1} \text{ (found in previous case)}; \text{ also,} \]
\[ x = \frac{DX}{XM} = \frac{DQ}{QA} \cdot \frac{AR}{RM} = q \cdot \frac{AR}{BR + BM} = \frac{q \cdot AR}{BR + \frac{1}{m+1} \cdot AB}, \text{ so } \]
\[ x = \frac{q}{\frac{BR}{AR} + \frac{1}{m+1} \cdot \frac{AB}{AR}} = \frac{q}{\frac{1}{r} + \frac{1}{m+1} \cdot \frac{r-1}{r}} = \frac{qr(m+1)}{r+m}. \text{ Then } \]
\[ y = \frac{CY}{YM} = \frac{CN}{NB} \cdot \frac{BR}{RM} = \frac{1}{n} \cdot \frac{RB}{RA} \cdot \frac{RA}{RM} = \frac{1}{rn} \cdot \frac{r(m+1)}{r+m} = \frac{m+1}{n(r+m)}, \]

\noindent
And the ratio $\frac{ME}{EP}$ takes the same value as before. \\

\noindent
$\bullet$ \textit{case 3.} $NQ \parallel AB$, then $r = 1$ by convention.
\medskip

\noindent
From parallelism, $QX \parallel AM \Rightarrow x = q$ and $YN \parallel MB \Rightarrow y = \frac{1}{n}$, therefore
\[ \frac{PE}{EM} = x \cdot \frac{p}{p+1} + y \cdot \frac{1}{p+1} = \frac{pq}{p+1} + \frac{1}{n(p+1)} = \frac{npq+1}{n(p+1)} = \frac{\gamma + m}{mn(p+1)}, \]
\[ \Rightarrow \frac{ME}{EP} = \frac{mn(p+1)}{\gamma + m} = \frac{mn(p+1)}{m+1} \cdot \frac{1+m}{\gamma + m} = \frac{mn(p+1)}{m+1} \cdot \frac{r+m}{\gamma r + m} \mid_{r=1}. \]
\end{proof}
\bigskip

Consider the initial configuration of \textit{Section 3}. One of our earlier observations was on the possible locations of points $F_{1}$, $F_{2}$ and $G_{1}$, $G_{2}$ on diagonals $AC$ and $BD$ respectively, depending on the value of $\gamma$. If we take into account their positions with respect to point $O$, there are several resulting configurations that may occur (see some of them below):
\clearpage

\begin{figure}[h]
\begin{center}
\includegraphics[width=0.85\linewidth]{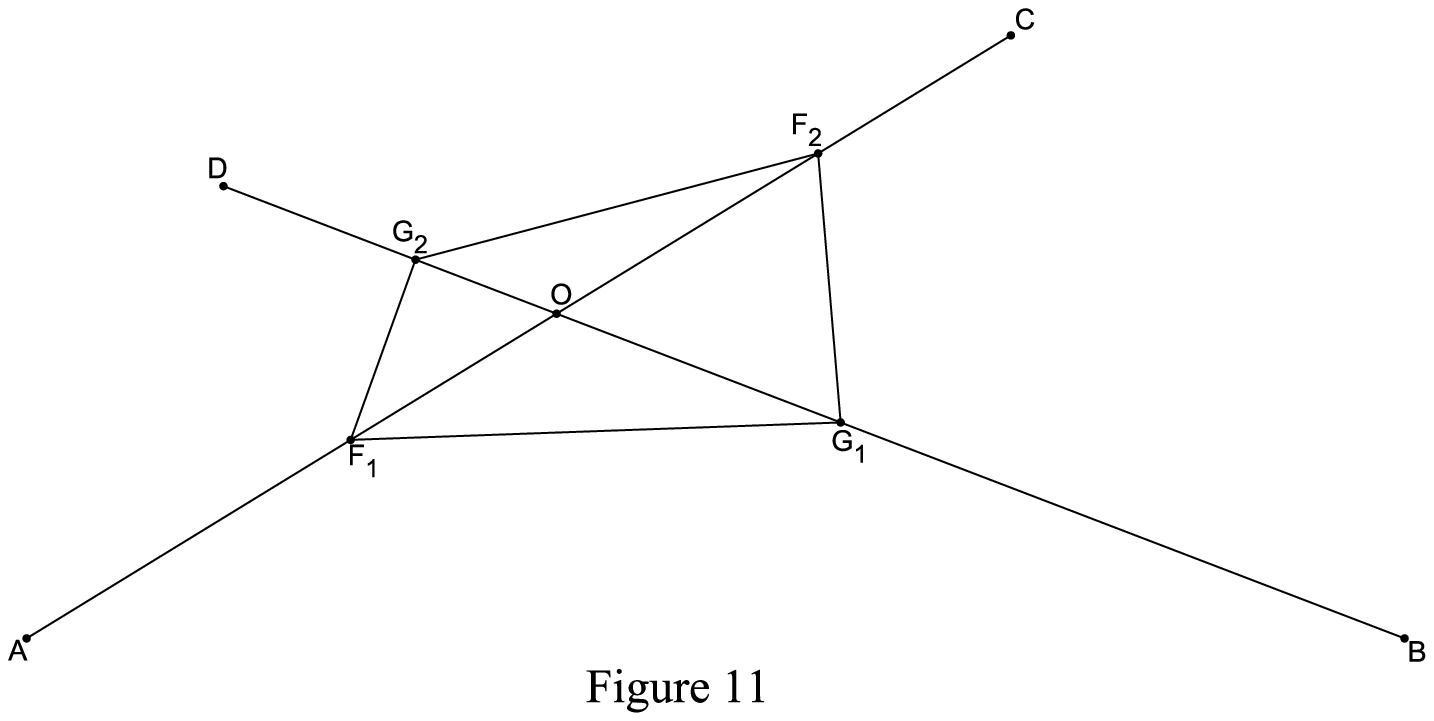}
\end{center}
\end{figure}
\vspace{-0.4in}

\begin{figure}[h]
\begin{center}
\includegraphics[width=0.85\linewidth]{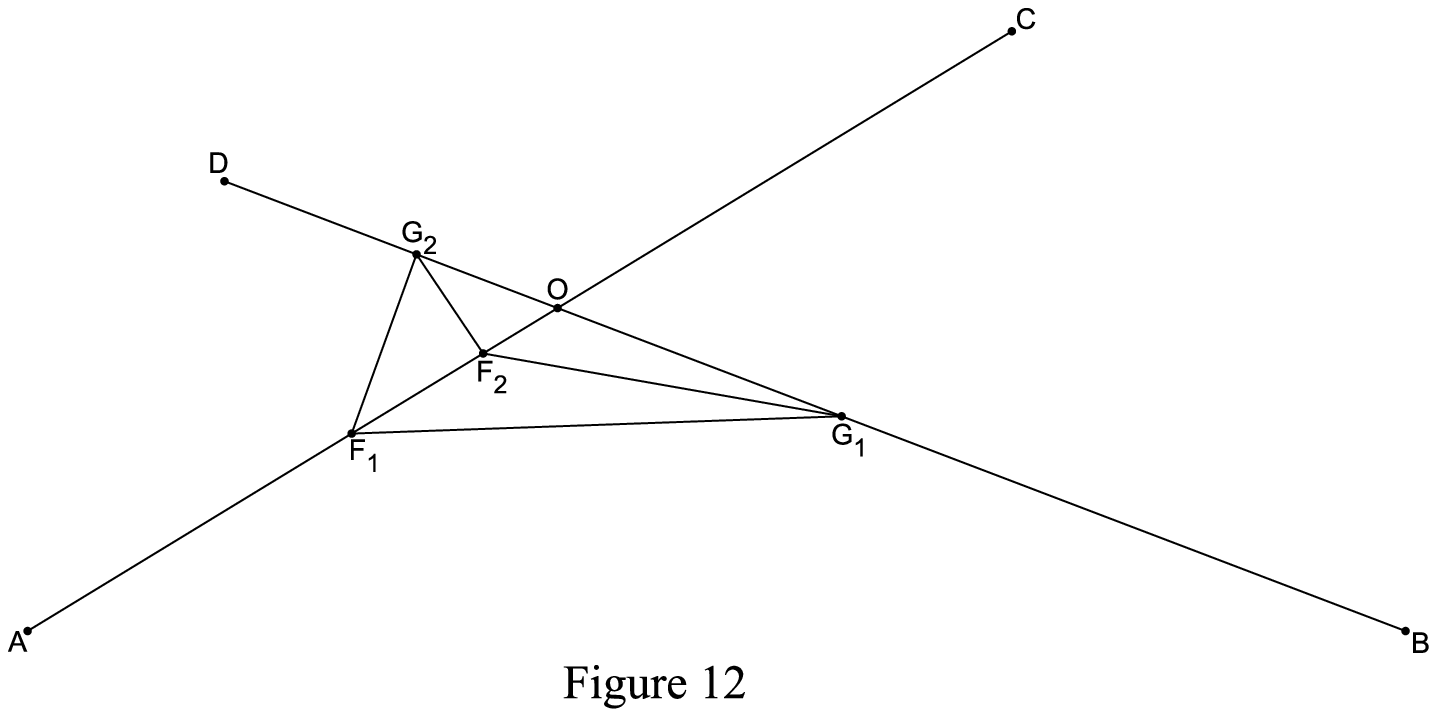}
\end{center}
\end{figure}

\noindent
\textbf{Question:} Can points $M_{1}$, $N_{1}$, $P_{1}$ and $Q_{1}$ be located as in the configuration from \textit{figure 13}? Can they describe anything else other than a convex quadrilateral? The answer is \textit{no}, and we will clearly see this in our next result.

\begin{figure}[htb]
\begin{center}
\includegraphics[width=\linewidth]{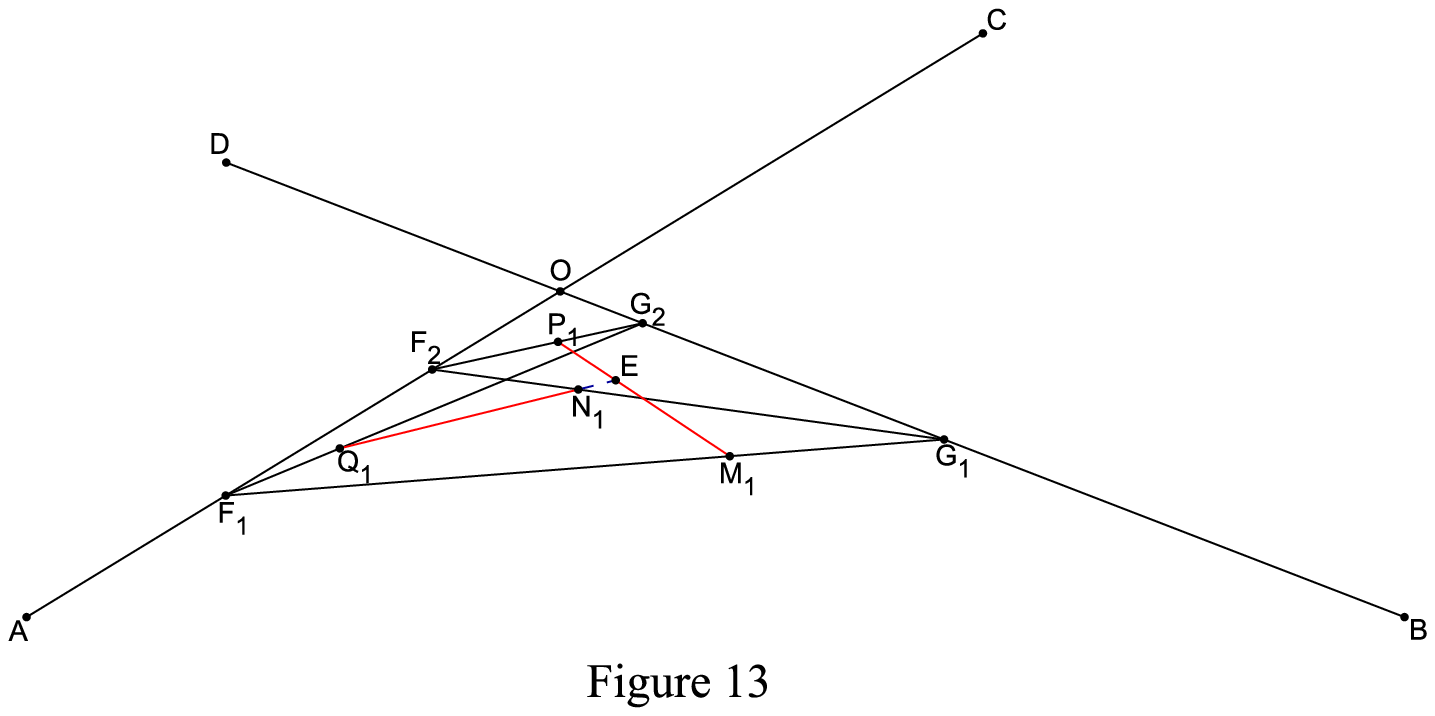}
\end{center}
\end{figure}

\begin{proposition}
If points $M_{1}$, $N_{1}$, $P_{1}$ and $Q_{1}$ are the intersection points from \textit{Theorem 11}, then $M_{1}N_{1}P_{1}Q_{1}$ is always a convex quadrilateral, possibly degenerate. 
\end{proposition}

\noindent
\begin{proof}
\, This statement is equivalent to saying that $M_{1}N_{1}P_{1}Q_{1}$'s diagonals meet at a point lying on both closed segments, i.e. $[M_{1}P_{1}] \cap [N_{1}Q_{1}] \neq \emptyset$. However $M_{1}, P_{1} \in MP$ and $N_{1}, Q_{1} \in NQ$, so $\{E\} = MP \cap NQ = M_{1}P_{1} \cap N_{1}Q_{1}$. Hence we only need to prove that $E \in [M_{1}P_{1}]$ and $E \in [N_{1}Q_{1}]$. By symmetry, proving the first part suffices.
\vspace{-0.1in}

\[ \text{We already found that } \frac{MM_{1}}{M_{1}P} = \frac{mn(p+1)}{m+1} \text{ (see \textit{Observation 3}), and also that} \]
\[ \frac{PP_{1}}{P_{1}M} = \frac{pq(m+1)}{p+1} \text{ (using \textit{Observation 1}), so } \frac{MP_{1}}{P_{1}P} = \frac{p+1}{pq(m+1)} = \frac{mn(p+1)}{m+1} \cdot \frac{1}{\gamma}. \]

\[ \text{Therefore, } E \in [M_{1}P_{1}] \Leftrightarrow \frac{ME}{EP} \in \left[ \min \left\{\frac{MM_{1}}{M_{1}P}, \frac{MP_{1}}{P_{1}P}\right\}, \max \left\{\frac{MM_{1}}{M_{1}P}, \frac{MP_{1}}{P_{1}P}\right\} \right]. \]

\noindent
Divide both sides by $\frac{MM_{1}}{M_{1}P} = \frac{mn(p+1)}{m+1}$ to get:

\[ E \in [M_{1}P_{1}] \Leftrightarrow \frac{r+m}{\gamma r + m} \in \left[ \min \left\{1, \frac{1}{\gamma}\right\}, \max \left\{1, \frac{1}{\gamma}\right\} \right]. \]
\smallskip

\noindent
$\bullet$ \textit{case 1.} $\gamma = 1$, then by \textit{Corollary 12}, $M_{1} \equiv P_{1} \equiv E$, and in particular $E \in [M_{1}P_{1}]$.
\smallskip

\noindent
$\bullet$ \textit{case 2.} $\gamma < 1$, implies $\gamma^2 r + \gamma m < \gamma r + \gamma m < \gamma r + m$ so $1 < \frac{r+m}{\gamma r + m} < \frac{1}{\gamma}$,
\[ \Rightarrow \frac{MM_{1}}{M_{1}P} < \frac{ME}{EP} < \frac{MP_{1}}{P_{1}P} \Rightarrow \text{the order is: } M-M_{1}-E-P_{1}-P \Rightarrow E \in [M_{1}P_{1}]. \]
\smallskip

\noindent
$\bullet$ \textit{case 3.} $\gamma > 1$, implies $\gamma^2 r + \gamma m > \gamma r + \gamma m > \gamma r + m$ so $1 > \frac{r+m}{\gamma r + m} > \frac{1}{\gamma}$,
\[ \Rightarrow \frac{MM_{1}}{M_{1}P} > \frac{ME}{EP} > \frac{MP_{1}}{P_{1}P} \Rightarrow \text{the order is: } M-P_{1}-E-M_{1}-P \Rightarrow E \in [M_{1}P_{1}]. \]
\smallskip

\noindent
$\therefore$ In conclusion, we have proved that $E \in [M_{1}P_{1}]$ and $E \in [N_{1}Q_{1}]$ follows by symmetry. Therefore, quadrilateral $M_{1}N_{1}P_{1}Q_{1}$ is convex.
\end{proof}
\bigskip

\begin{remark}
So far we assumed throughout the article that $ABCD$ is convex. However both \textit{Theorem 4} and \textit{Theorem 11} still hold when we change the quadrilateral to be concave, or even complex (self-intersecting). These two situations are left to the reader and should be treated separately, while keeping in mind that the auxiliary lemmas and the proofs remain the same.
\end{remark}
\medskip

\begin{remark}
If $ABCD$ is either concave or self-intersecting, then some of the results derived in both \textit{Section 2} and \textit{Section 3} do not hold any longer. For example, we may notice in the configuration in \textit{figure 14} that the quadrilateral $M_{1}N_{1}P_{1}Q_{1}$ is self-intersecting and not convex, as suggested by \textit{Proposition 14}.
\end{remark}

\begin{figure}[htb]
\begin{center}
\includegraphics[width=\linewidth]{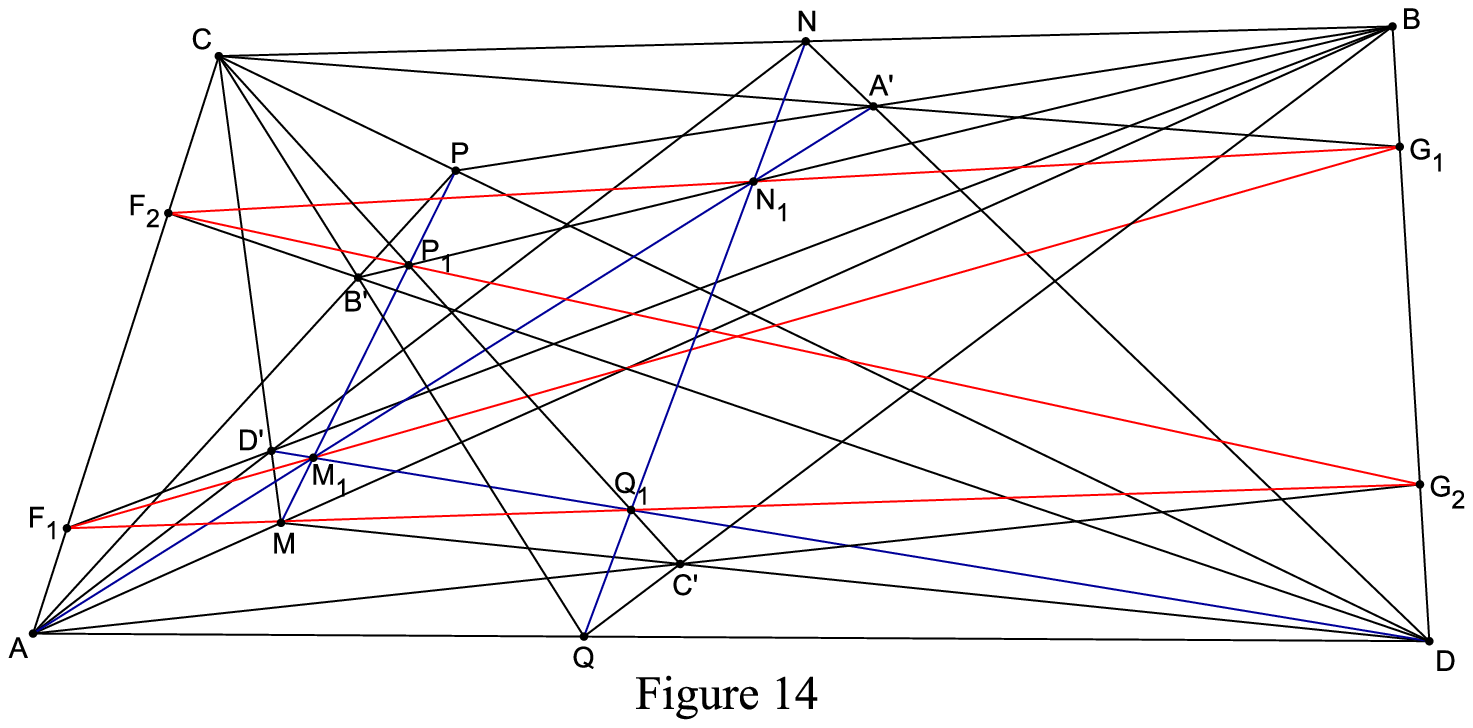}
\end{center}
\end{figure}
\bigskip

 
\medskip

Andrei Cozma: $25$ Grandale Street, Manchester M$145$WQ, United Kingdom

\textit{E-mail address}: \verb"eagle_nest23@yahoo.com"

\end{document}